\documentclass{amsart}
\usepackage[latin1]{inputenc}
\usepackage{amssymb}
\newtheorem{thm}{Theorem}

\newtheorem{prop}[thm]{Proposition}

\theoremstyle{definition}

 \numberwithin{equation}{section}

\newcommand{\To}{\longrightarrow}
\begin{document}

\title[]{Automorphisms in spaces of continuous functions on Valdivia compacta}%
\author{Antonio Avilés and Yolanda Moreno}

 \address{University of Paris 7, Equipe de
Logique Mathématique, UFR de Mathématiques, 2 place Jussieu, 75251
Paris, France}
\address{Departamento de Matemáticas, Universidad de Extremadura, Escuela Politécnica, Av. Universidad, s/n. 10.071 - Cáceres}%

\email{avileslo@um.es, ymoreno@unex.es}

\thanks{The first author was supported
by a Marie Curie Intra-European Felloship MCEIF-CT2006-038768
(E.U.) and research projects MTM2005-08379 and S\'{e}neca
00690/PI/04 (Spain)}

\begin{abstract}
We show that there are no automorphic Banach spaces of the form
$C(K)$ with $K$ continuous image of Valdivia compact except the
spaces $c_0(\Gamma)$. Nevertheless, when $K$ is an Eberlein
compact of finite height such that $C(K)$ is not isomorphic to
$c_0(\Gamma)$, all isomorphism between subspaces of $C(K)$ of size
less than $\aleph_\omega$ extend to
automorphisms of $C(K)$.
\end{abstract}

\maketitle

\section*{Introduction}

A Banach space $X$ is said to be automorphic if for every
isomorphism $T:Y_1\To Y_2$ between two (closed) subspaces of $X$
with $dens(X/Y_1)=dens(X/Y_2)$ there exists an automorphism
$\tilde{T}:X\To X$ which extends $T$, that is,
$\tilde{T}|_{Y_1}=T$. It has been shown in \cite{MorPli} that a
necessary condition for a Banach space $X$ to be automorphic is to
be extensible, which means that for every subspace $E\subset X$
and every operator $T:E\To X$, there exists an operator
$\tilde{T}:X\To X$ that extends $T$. Clearly every Hilbert space
$\ell_2(\Gamma)$ is automorphic and on the other hand,
Lindenstrauss and Rosenthal \cite{c0automorphic} have proven that
$c_0$ is automorphic and also that $\ell_\infty$ has a partial
automorphic character, namely that isomorphisms $T:Y_1\To Y_2$ can
be extended provided that $\ell_\infty/Y_i$ is nonreflexive for
$i=1,2$, though $\ell_\infty$ is not automorphic. Moreno and
Plichko \cite{MorPli} have recently shown that $c_0(\Gamma)$ is
automorphic for every set $\Gamma$. It remains open the question
posed in \cite{c0automorphic} whether the only automorphic
separable Banach spaces are $\ell_2$ and $c_0$ and also the more
general question whether all automorphic Banach spaces are
isomorphic either to $\ell_2(\Gamma)$ or to $c_0(\Gamma)$ for some
set $\Gamma$. Our aim in this note is to address this latter
problem for the case of Banach spaces $C(K)$ of continuous
functions on compact spaces. Must be an automorphic $C(K)$ space
isomorphic to $c_0(\Gamma)$? We provide a positive answer to this
problem in the case when $K$ is a continuous image of a Valdivia
compact, which is a large class of compact spaces originated from
functional analyisis and which includes for example all Eberlein
and all dyadic compact
spaces.\\

Namely, a compact space is said to be a Valdivia compact if it is
homeomorphic to some $K\subset \mathbb{R}^\Gamma$ in such a way
that the elements of $K$ of countable support are dense in $K$
(the support of $x\in \mathbb{R}^\Gamma$ is the set of nonzero
coordinates). If such $K$ can be found so that all elements of $K$
have countable support, then the compact is said to be a Corson
compact, and if moreover it can be taken so that $K\subset
c_0(\Gamma)\subset \mathbb{R}^\Gamma$, then it is called an
Eberlein compact.\\

\begin{thm}\label{Valdiviaextensible}
Let $K$ be a continuous image of a Valdivia compact. If $C(K)$ is
extensible, then $C(K)$ is isomorphic to $c_0(\Gamma)$.\\
\end{thm}

Although it is a standard notion, we recall now what a scattered
compact is. The derived space of a topological space $X$ is the
space $X'$ obtained by deleting from $X$ its isolated points. The
derived sets $X^{(\alpha)}$ are defined recursively setting
$X^{(0)}=X$, $X^{(\alpha+1)}=[X^{(\alpha)}]'$ and
$X^{(\beta)}=\bigcap_{\gamma<\beta}X^{(\gamma)}$ for $\beta$ a
limit ordinal. The space $X$ is scattered if
$X^{(\alpha)}=\emptyset$ for some $\alpha$, and in this case the
minimal such $\alpha$ is called the height of $X$.\\

In Theorem~\ref{Valdiviaextensible}, if $K$ is not a scattered
compact of finite height, then the extensible property fails at
the separable level, meaning that in that case $C(K)$ contain both
a complemented and non complemented copy of the same separable
space. The most delicate case of Theorem~\ref{Valdiviaextensible}
happens when $K$ is a scattered compact of finite height and it
relies on some results of~\cite{ArgyrosandSpain}. In this
situation $K$ is indeed an Eberlein compact and the special
behavior of $c_0(\Gamma)$ when $|\Gamma|<\aleph_\omega$ studied in
\cite{ArgyrosandSpain}, \cite{BelMarK} and \cite{GodKalLan}
combined with the general results about $c_0(\Gamma)$ from
\cite{MorPli} yield that we have the extensible and automorphic
properties for subspaces of density less than $\aleph_\omega$.\\

\begin{thm}\label{almostautomorphic}
Let $K$ be an Eberlein compact of finite height.
\begin{enumerate}
\item For every isomorphism $T:Y_1\To Y_2$ between two subspaces
of $X$ with $dens(Y_1)=dens(Y_2)<\aleph_\omega$ and
$dens(C(K)/Y_1) = dens(C(K)/Y_2)$ there exists an
automorphism $\tilde{T}:C(K)\To C(K)$ that extends $T$.\\

\item For every subspace $Y\subset C(K)$ with
$dens(Y)<\aleph_\omega$ and every operator $T:Y\To C(K)$, there
exists an operator $\tilde{T}:C(K)\To C(K)$ that extends $T$.
\end{enumerate}
\end{thm}

Only recently, Bell and Marciszewski~\cite{BelMarK} have
constructed an Eberlein compact space of height 3 and weight
$\aleph_\omega$ that is not isomorphic to $c_0(\Gamma)$, where
$|\Gamma|=\aleph_\omega$. It was shown by Godefroy, Kalton and
Lancien~\cite{GodKalLan} that if $K$ is an Eberlein compact of
finite height and weight less than $\aleph_\omega$,
then $C(K)$ is isomorphic to $c_0(\Gamma)$, cf. also \cite{BelMarK} and \cite{Marciszewskic0}.\\

The most typical example of scattered compact which is not
Eberlein is a Mr\'{o}wka space, that is, a separable uncountable
scattered compact space $K$ of height three and $|K^{(2)}|=1$. In
this case $C(K)$ is not extensible, cf. Proposition~\ref{Mrowka}.
However, it is unclear to us whether there may exist a scattered
compact space such that $C(K)$ is extensible but not isomorphic to
any
$c_0(\Gamma)$.\\

This research was done when both authors were visiting the
National Technical University of Athens. We are specially indebted
to Spiros Argyros for his valuable help and suggestions.\\

\section*{Proof of Theorem~\ref{Valdiviaextensible}}
Let us first observe that if $C(K)$ is an extensible space then $K$
does not contain any copy of the ordinal interval
$[0,\omega^\omega]$. Indeed, if we had $[0,\omega^\omega]\subset K$,
by the Borsuk-Dugundji extension theorem, $C[0,\omega^\omega]$ is a
complemented subspace of $C(K)$. But it is known (see \cite{Pelczy})
that $C[0,\omega^\omega]$ contains an uncomplemented copy of itself,
so $C(K)$ contains both complemented and uncomplemented copies of
$C[0,\omega^\omega]$ and so $C(K)$ is not extensible.\\

A result of Kalenda asserts that a continuous image of a Valdivia
compact which does not contain the ordinal interval $[0,\omega_1]$
is a Corson compact. Any Corson compact is monolithic (that is,
every separable closed subset has countable weight) and for
monolithic spaces we have the following, which is probably a known
fact:\\

\begin{prop}
Let $K$ be a monolithic compact which does not contain any copy of
$[0,\omega^\omega]$, then $K$ is a scattered compact of finite
height.\\
\end{prop}

Proof: First, if $K$ was not scattered, then there is a continuous
surjection from $K$ onto the unit interval, $f:K\To [0,1]$. For
every rational point $q\in [0,1]$ we choose $x_q\in K$ with
$f(x_q)=q$. Let $L$ be the closure of the set of points $\{x_q :
q\in \mathbb{Q}\cap [0,1]\}$, which is a metrizable compact since
$K$ is monolithic. Moreover, $f$ maps $L$ onto $[0,1]$, so $L$
contains a perfect compact set, which being metrizable, contains a
further copy of the Cantor set $\{0,1\}^\mathbb{N}$ and in
particular a copy of
$[0,\omega^\omega]$, contrary to our assumption. Thus, we proved that $K$ must be scattered.\\

Suppose now that $K$ was a scattered compact of infinite height.
For every $n\in\mathbb{N}$ let $A_n = K^{(n)}\setminus K^{(n+1)}$
be the $n-th$ level of $K$. Since the height of $K$ is infinite,
$A_n$ is nonempty for every $n\in\mathbb{N}$. Indeed, $A_n$ is an
infinite set which is dense in $K^{(n)}$. We observe that for
$n>0$, every element $x\in A_n$ is the limit of a sequence of
elements of $A_{n-1}$ (take $U$ a clopen set which isolates $x$
inside $K^{(n)}$, then $U\cap A_{n-1}$ is infinite and any
sequence contained in $U\cap A_{n-1}$ converges to $x$). For every
$n$ we take countable sets $B_{n,n}\subset A_n$, $B_{n,n-1}\subset
A_{n-1}$, $\ldots$ $B_{n,0}\subset A_0$ in the
following way:\\

\begin{itemize}
\item $B_{n,n}$ is an arbitrary countably infinite subset of
$A_n$.\\

\item $B_{n,n-1}$  is a countable subset of $A_{n-1}$ such that
every element of $B_{n,n}$ is the limit of a sequence of elements
of $B_{n,n-1}$.\\

\item $B_{n,k}$ is a countable subset of $A_{k}$ such that every
element of $B_{n,k+1}$ is the limit of a
sequence of elements of $B_{n,k}$.\\
\end{itemize}

Let $B_n = B_{n,n}\cup B_{n,n-1}\cup\cdots\cup B_{n,0}$. Notice
that $B_n$ is a scattered topological space of height $n+1$ with
$B_n^{(k)} = B_{n,n}\cup B_{n,n-1}\cup\cdots\cup B_{n,k}$. Let
$L=\overline{\bigcup_{n=0}^\infty B_n}$. The compact $L$ is a
scattered compact of infinite height and it is moreover metrizable
because $K$ is monolithic. Any metrizable scattered compact is
homeomorphic to an ordinal interval, and since $L$ has infinite
height, $[0,\omega^\omega]\subset L$.$\qed$\\

In order to prove Theorem~\ref{Valdiviaextensible} we shall assume
by contradiction that there exists some compact space $K$ which is
a continuous image of Valdivia compact, with $C(K)$ extensible and
not isomorphic to $c_0(\Gamma)$. The previous discussion shows
that any such $K$ must be scattered compact of finite height.
Hence, we can choose one such compact $K_0$ of minimal height. We
shall work with this
$K_0$ towards getting a contradiction.\\

Let $\Delta$ be the set of isolated points of $K_0$, so that
$K'_0=K_0\setminus\Delta$, $K_0=K'_0\cup\Delta$. We consider the
restriction operator $S:C(K_0)\To C(K_0')$ for which $ker(S) =
c_0(\Delta)$ and we have a short exact sequence

$$(\star)\ \ 0\To c_0(\Delta)\To C(K_0)\To C(K'_0)\To 0.$$

By~\cite[Theorem~1.2]{ArgyrosandSpain}, there exists
$\tilde{\Delta}\subset \Delta$ with $|\tilde{\Delta}|=|\Delta|$
such that $c_0(\tilde{\Delta})$ is complemented in $C(K_0)$. Since
$C(K_0)$ is extensible, it follows that, being
$c_0(\tilde{\Delta})$ a complemented subspace of $X$, also
$c_0(\Delta)$ is a complemented subspace of $C(K_0)$. Therefore,
the short exact sequence $(\star)$ splits and we have

$$C(K_0) = c_0(\Delta)\oplus C(K'_0).$$

In particular, $C(K_0')$ is a complemented subspace of $C(K_0)$
and therefore $C(K_0')$ is also extensible. Moreover, $K_0'$ has
height one unit less than the height of $K_0$, so by the
minimality property used to choose $K_0$ we conclude that
$C(K_0')$ is isomorphic to $c_0(\Gamma)$ for some $\Gamma$. But
then

$$C(K_0)\cong c_0(\Delta)\oplus C(K_0')\cong
c_0(\Delta)\oplus c_0(\Gamma)\cong c_0(\Delta\cup\Gamma),$$

a contradiction since $C(K_0)$ was not isomorphic to any
$c_0(\Lambda$). This finishes the proof of
Theorem~\ref{Valdiviaextensible}.\\

Let us note that we did not use the full strength of the
assumption of $C(K)$ being extensible in the hypothesis of
Theorem~\ref{Valdiviaextensible}. We only needed that $C(K)$ does
not contain both complemented and uncomplemented copies of the
same space $X$, for the spaces $X=C[0,\omega^\omega]$ and
$X=c_0(\Gamma)$.\\

We include now the proof that Mrowka compacta do not provide
extensible Banach spaces, which uses similar ideas as in the
preceding arguments:\\

\begin{prop}\label{Mrowka}
Let $K$ be a Mrowka space. Then $C(K)$ contain both complemented
and uncomplemented copies of $c_0$.\\
\end{prop}

Proof: $K$ contains convergent sequences, that is, a copy of
$[0,\omega]$, so by the Borsuk-Dugundji extension theorem, it
contains a complemented copy of $C[0,\omega]\cong c_0$. On the
other hand, let $\Delta$ be the countable set of the isolated
points. Like above, $c_0(\Delta)$ is the kernel of the restriction
operator $C(K)\To C(K')$. It is well known that $c_0(\Delta)$ is
not complemented in $C(K)$ in this case. One argument to see this
is the following: Suppose $c_0(\Delta)$ was complemented in
$C(K)$. Then $C(K)\cong c_0(\Delta)\oplus C(K')$. The space $K'$
is homeomorphic to the one point compactification of a discrete
set $\Gamma$, so $C(K')\cong c_0(\Gamma)$ and $C(K)\cong
c_0(\Delta)\oplus c_0(\Gamma)$. This implies $C(K)$ is a weakly
compactly generated space, and therefore $K$ is an Eberlein
compact. Every separable Eberlein compact has countable weight and
a Mrowka space is separable but has uncountable weight (we refer
to \cite{FabianWA} for
reference to standard properties of weakly compactly generated spaces and Eberlein compact spaces).$\qed$\\


\section*{Proof of Theorem~\ref{almostautomorphic}}

Let us first note that a continuous image $K$ of Valdivia compact
which is scattered compact of finite height is an Eberlein
compact. We use again Kalenda's result \cite{Kalenda} that $K$
must be either Corson or contain a copy of $[0,\omega_1]$, and the
latter possibility cannot happen since $K$ has finite height. It
is a result of Alster~\cite{Alster} that every scattered Corson
compact is an
Eberlein compact.\\

We state now a result from~\cite{GodKalLan} mentioned in the
introduction:\\

\begin{thm}[Godefroy, Kalton, Lancien]\label{Ksmallweight}
If $Q$ is an Eberlein compact of finite height and
$w(Q)=\aleph_m<\aleph_\omega$, then $C(Q)$ is isomorphic to
$c_0(\aleph_m)$.\\
\end{thm}

In the view of this and of the fact that $c_0(\Gamma)$ is an
automorphic and hence also extensible space, we are concerned in
Theorem~\ref{almostautomorphic} with the case when $w(K)\geq
\aleph_\omega$. So let $K$ be an Eberlein compact of finite height
and weight not lower than $\aleph_\omega$ and let $T:Y_1\To Y_2$
be an isomorphism between subspaces of $C(K)$ such that
$dens(Y_1)=dens(Y_2)=\aleph_n<\aleph_\omega$. Our aim is to find
an automorphism of $C(K)$ that extends $T$.\\

 Let $Z$ be a subspace of $C(K)$ of density character $\aleph_{n+1}$
 such that $Y_1+Y_2\subset Z$. We define an equivalence relation $\sim$ on $K$ in the following
 way:

 $$p\sim q \iff y(p)=y(q) \text{ for every }y\in Z.$$

The quotient $L=K/\sim$ with the quotient topology is a compact
space and the quotient map $K\To L$ a continuous surjection which
allows us to view $C(L)$ as a subspace of $C(K)$ such that
$Z\subset C(L)$. Moreover, since the space $Z$ separates the
points of $L$ and has density character $\aleph_{n+1}$,
$w(L)=\aleph_{n+1}$. Now, by Theorem~\ref{Ksmallweight}, $C(L)$ is
isomorphic to $c_0(\aleph_{n+1})$ and we know by~\cite{MorPli}
that this space is automorphic. Hence, since
$dens(C(L)/Y_1)=\aleph_{n+1}=dens(C(L)/Y_2)$, there exists an
automorphism $\hat{T}:C(L)\To C(L)$ that extends $T$. Finally, by
\cite[Theorem 1.1]{ArgyrosandSpain} every copy of
$c_0(\aleph_{n+1})$ in a weakly compactly generated space is
complemented, so $C(L)$ is complemented in $C(K)$ and this allows
us to obtain an automorphism $\tilde{T}:C(K)\To C(K)$ that extends
$\hat{T}$. This finishes the proof of part (1) of
Theorem~\ref{almostautomorphic}.\\

Part (2) of Theorem~\ref{almostautomorphic} is a consequence of
part (1) by \cite[Theorem 3.1]{MorPli} (this theorem only states
that an automorphic space must be extensible but the proof shows
that if the automorphic property holds for a given subspace then
so does the extensible property). Alternatively, part (2) can also
be proven by an argument which is completely analogous to that of
part (1).


\begin{thebibliography}{10}

\bibitem{Alster} K. Alster \emph{Some remarks on Eberlein compacts}
Fundam. Math. 104 (1979) 43-46.

\bibitem{ArgyrosandSpain} S. A. Argyros, J. F. Castillo, A. S.
Granero, M. Jiménez and J. P. Moreno,\emph{Complementation and
embeddings of $c_0(I)$ in Banach spaces}. Proc. London Math. Soc.
(3) 85 (2002) 742-768.

\bibitem{BelMarK}
M. Bell and W. Marciszewski, \emph{On scattered Eberlein compact
spaces}. Israel J. Math. 158 (2007), 217-224.

\bibitem{FabianWA}
M.~Fabian, \emph{G\^ateaux differentiability of convex functions
and topology},
  Canadian Mathematical Society Series of Monographs and Advanced Texts, John
  Wiley \& Sons Inc., New York, 1997, Weak Asplund spaces, A Wiley-Interscience
  Publication.

\bibitem{GodKalLan} G. Godefroy, N. Kalton and G. Lancien. \emph{Subspaces of
$c_0(\mathbb{N})$ and Lipschitz isomorphisms}. Geom. Funct. Anal.
10, No.4, (2000) 798-820.

\bibitem{Kalenda} O. F. K. Kalenda, \emph{Valdivia compact spaces in topology and
Banach space theory}. Extr. Math. 15, No.1, (2000) 1-85.

\bibitem{c0automorphic}
J. Lindenstrauss and H. P. Rosenthal, \emph{Automorphisms in
}$c_0$, $l_1$ \emph{and} $m$, Israel J. Math. 7 (1969), 227-239.

\bibitem{Marciszewskic0} W. Marciszewski, \emph{On Banach spaces
$C(K)$ isomorphic to $c_0(\Gamma)$}, Studia Math. 156 (2003),
295-302.

\bibitem{MorPli}
Y. Moreno and A. Plichko, \emph{On extending automorphisms in
Banach spaces}. To appear in Israel J. Math.

\bibitem{Pelczy}
A. Pelczynski, \emph{Linear extensions, linear averagings, and their
applications to linear topological classification of spaces of
continuous functions}, Dissertationes Math. 58, 1968.

\end{thebibliography}
\end{document}